\documentclass[12pt]{article}

\usepackage{texdraw}
\usepackage{graphicx}
\usepackage{amssymb}

\newcommand{\be}{\begin{equation}}
      \newcommand{\ee}{\end{equation}}
      \newcommand{\ba}{\begin{eqnarray}}
       \newcommand{\ea}{\end{eqnarray}}
\newcommand{\ban}{\begin{eqnarray*}}
\newcommand{\ean}{\end{eqnarray*}}

\newcommand{\RR}{\mathbb{R}}
\newcommand{\EE}{\mathbb{E}}

\newcommand{\ZZ}{\mathbb{Z}}
\newcommand{\NN}{\mathbb{N}}

\newcommand{\sect}[1]{\section{#1} \setcounter{equation}{0}}

\begin{document}
\title{The Topology of Open Manifolds with Nonnegative Ricci Curvature}

\author{Zhongmin Shen and Christina Sormani}

\date{}
\maketitle

\noindent{\bf Abstract:} We survey all results
concerning the topology of open
manifolds with $Ricci \ge 0$ that have
no additional conditions other than 
restrictions to the dimension, volume growth or diameter
growth of the manifold.  We will also present relevant
examples and list open problems.

\sect{Introduction}

Throughout this survey article, $M^n$ is a
complete noncompact
$n$ dimensional manifold 
with $Ricci \ge 0$ where $n\ge 3$.
Note that Cohn-Vossen proved
$M^2$ is either diffeomorphic
to $\RR^2$ or flat \cite{Cv}.  
DeTurck proved that locally
one can prescribe Ricci curvature, so all obstructions to the
existence of such a metric are topological \cite{Dt}.
Lohkamp has demonstrated that there are no topological 
obstructions for manifolds
with negative Ricci curvature \cite{Lok} . 

In Section~\ref{geods} we provide background on
geodesics, Ricci curvature and the Cheeger-Gromoll 
Splitting Theorem.  
We then describe warped products and the 
examples of Nabonnand, Wei and Wilking.
This leads into the Loops to Infinity theorem of the second author.
Section~\ref{geods} closes with a complete
classification of the codimension one homology, $H_{n-1}(M,Z)$, 
and a small restriction to $H_{n-2}(M,Z)$
derived by the two authors in 2000.

In Section~\ref{volume} we  
present the Bishop-Gromov Volume Comparison Theorem and introduce
Milnor's Conjecture that the $\pi_1(M)$
is finitely generated.  We describe results of Milnor,
Gromov and Wilking in this direction.
Partial solutions of the Milnor Conjecture
requiring additional hypothesis on volume or diameter by
Li, Anderson and the second author are presented in detail
as well as the Perelman Contractibility Theorem.  We also
mention the Cheeger-Colding Diffeomorphism Theorem.

Section~\ref{examples} begins with thorough descriptions of
examples of $M^n$ with infinite topological type by Sha-Yang.
Anderson-Kronheimer-LeBrun, and Menguy.  We also review
Wraith's surgery techniques and Perelman's building blocks.
We close Section~\ref{examples}
with a description of  results
of Nash, Berard-Bergery, Otsu, Anderson, and Belegradek-Wei
exploring which vector bundles admit metrics
with $Ricci >0$ and $Ricci \ge 0$.

In Section~\ref{3D} we present Schoen-Yau's proof that
three manifolds with $Ricci >0$ are diffeomorphic to $\RR^3$.  
The classification of the topology of $M^3$ with only $Ricci \ge 0$
is an open problem.   We describe partial results by Schoen-Yau,
Shi, Zhu, Meeks-Simon-Yau, 
Anonov-Burago-Zalgaller, Anderson-Rodriguez and Zhu.  
Zhu in fact has shown $M^3$ is contractible as long as the volume
grows like $r^3$.  However, Milnor's Conjecture remains
open even in dimension three!
Section~\ref{3D} closes with a description of a potential Milnor
Counter Example: the dyadic solenoid complement.  

Section~\ref{open} describes further open problems.  It includes a list 
of the qualitative properties of $M^n$ for those wishing to search
for new examples.  
We close the paper with thanks to the many geometers 
and topologists who assisted us.

There are many beautiful results on the topology of open manifolds with
nonnegative Ricci curvature which have additional conditions on either
the Busemann function, injectivity radius, conjugacy radius or some other
geometric constraint.  However, we were unable to include these results here.
We have also had to leave out the related theory of compact manifolds
including many relevant examples.  To keep the bibliography shorter
than five pages we only refer to the primary articles we are surveying
and not the important papers cited within those articles.  
We hope that this survey will prove useful to everyone interested
in entering this area rich in open problems.

\sect{Geodesics and $H_{n-1}(M,Z)$} \label{geods}

In this section we provide some intuitive understanding of
geodesics.  Throughout, geodesics will be parametrized by
arclength.  We begin without the assumption of Ricci curvature.

A {\bf ray} is a geodesic $\gamma:[0,\infty) \longmapsto M$ such that
$$
d(\gamma(t), \gamma(s))=|s-t|\qquad 
\forall s,t\in[0,\infty)
$$
Recall that a geodesic fails to minimize after its
first cut or conjugate point, so as soon as there
is more than one path between a pair of points on
a geodesic it is no longer a ray.
On a paraboloid, $z=x^2+y^2$, any geodesic running radially
outward from the basepoint $(0,0,0)$ is a ray.  
In fact {\em every complete 
noncompact Riemannian 
manifold contains a ray}.

A {\bf line} is a geodesic 
$\gamma:(-\infty, \infty)\longmapsto M^n$, such that:
$$
d(\gamma(t),\gamma(s))=|t-s| 
\forall \,\, t,\, s\in (-\infty,\, \infty).
$$
Paraboloids have no lines while cylinders have a collection
of parallel lines in their so-called ``split'' direction.

We say a manifold had $k$ {\bf ends} if for every 
sufficiently large compact set, $K$,
$M\setminus K$ has $k$ unbounded components.
So a paraboloid has one end, a cylinder has two ends
and a Riemann surface with three punctures or cusps has
three ends.
{\em A complete manifold with two or more ends contains a line.}

Ricci curvature, on the other hand, is a locally defined concept:
Let $p\in M^n$,  $v \in TM_p$ with $|v|=1$,
{\bf{Ricci curvature}} is defined
\be \label{sum}
Ric_p(v,v)=\sum_{i=1}^{n-1}<R(e_i,v)v,e_i>
\ee
where $\,\,v, e_1, e_2,...e_{n-1}\,\,$ are 
orthonormal and $R$ is the 
sectional curvature tensor.
We say $M$ has $Ricci \ge 0$ if
\be \label{nonneg}
Ric_p(v,v)\ge 0 \,\,\forall p\in M\,\,\forall v\in TM_p
\ee
and it has $Ricci >0$ if $\ge$ is replaced by $>$ in (\ref{nonneg}).

\begin{figure}[htbp]
\vspace{.2cm}
\begin{center}
\includegraphics[height=1in ]{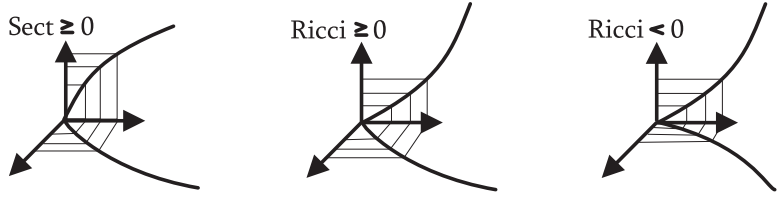}
\end{center}
\caption{The bending of geodesics under curvature bounds.} 
\label{bendapart}
\end{figure}

Intuitively, the sectional curvature measures how much 
geodesics bend together.  Thus if $(R(v,w)w,v)\ge 0$ then
the geodesics $C_v(t)=exp_p(tv)$ and $C_w(t)=exp_p(tw)$
starting at $p$ in the directions $v$ and $w$ respectively,
tend to bend towards each other or diverge at most linearly.
On the other hand if only $Ricci_p(v,v)\ge 0$ then 
by (\ref{sum}) a pair of geodesics may bend apart
as long as other 
geodesics bend together.  See Figure~\ref{bendapart}.

{\bf Cheeger-Gromoll Splitting Theorem (1971):}
{\em If $M^n$ contains a line, 
then $M^n$ splits isometrically:
$$
M^n= (-\infty, \infty) \times N^{n-k}
$$
with the metric, $g_M=dr^2 + g_N$.
} \cite{ChGl2}

This theorem can be intuitively understood
as saying all geodesics must remain parallel to the line 
because if some were to bend away 
then others would have to bend inward 
causing a cut point.
The actual proof of this theorem uses an arguement
involving subharmonic and superharmonic functions (c.f. \cite{EscHe}).

{\em Note that Cheeger-Gromoll's Splitting theorem
implies that $M^n$ has at most two ends and $M^n$
with $Ricci >0$ has only one end.}

A {\bf warped product}  $M= [0,\infty) \times_f S^1$
is a manifold with the metric 
$$
g_M=dr^2 + f^2(r)d\theta^2.
$$
It is a smooth manifold with one end
if $f(r)>0$ for $r>0$,
$f(0)=0$, $f'(0)=1$, and $f''(0)=0$.

If $f(r)=r$ this is Euclidean space.
If $f(r)=sin(r)$ this is a sphere and
it closes up at $r=\pi$.
In fact, $[0,\infty)\times_f S^k$ is a manifold 
with positive sectional curvature
if $f''(r) <0$.
Curves of the form $\gamma(t)=(t, \theta_0)$
are geodesics in this warped product.
One can see that these geodesics are bending together
when $f''(r)>0$.  
In order to construct more interesting
examples with positive Ricci curvature that have some
negative sectional curvature, one needs to warp 
the manifold with more than one function.

A {\bf doubly warped product}
$
[0,\infty)\times_h S^2 \times_f S^1
$
has  the metric
\be
dr^2+h^2(r)g_{S^2} + f^2(r) d\theta^2.
\ee
Its radial Ricci curvature is
\be  \label{eqnwarp}
Ric(\frac{\partial }{\partial r},\frac{\partial }{\partial r})
=-\frac{f''(r)}{f(r)} -2\frac{h''(r)}{h(r)}.
\ee
Notice now how some curves now may bend apart (say with
$f''(r)>0$ as long as others bend together (with $h''(r)<0$).

\begin{figure}[htbp]
\begin{center}
\includegraphics[height=1.5in ]{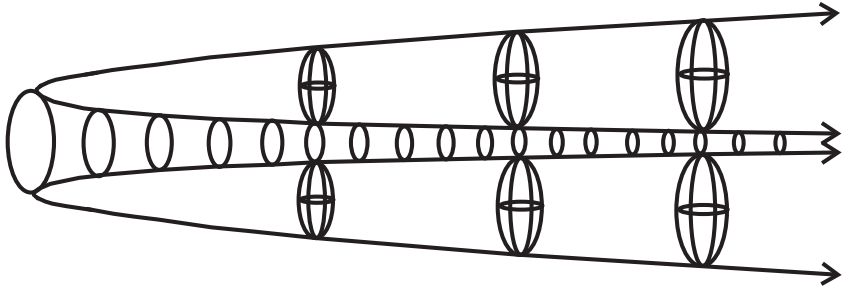}
\end{center}
\caption{Nabonnand's Example} \label{nab1}
\end{figure}

{\bf Example of Nabonnand (1980):} 
{\em There is a doubly warped product
\be
M^4=[0,\infty)\times_h S^2 \times_f S^1
\ee
with $Ricci>0$, $h(0)=0$, $h'(0)=1$, and $h''(0)=0$
but $f(0) > 0$.  In particular $M^4$ is 
diffeomorphic to $\RR^3\times S^1$ and has
fundamental group $\pi_1(M^4)=\ZZ$.
See Figure~\ref{nab1}.} \cite{Nab}

Note that in Nabonnand's example
the loop at $r=0$ is not contractible.
However, it is homotopic to shorter and 
shorter loops diverging to infinity.

{\bf Examples of Wei (1988):}
{\em For any discrete nilpotent group, $G$,
there is an $M^n$ with fundamental group, $\pi_1(M^n)=G$:
$$
M=[0,\infty)\times_h S^k \times_f N
$$
with $h(0)=0$ but $f(0)\neq 0$,
positive Ricci curvature, and
fundamental group $\pi_1(N)=G$.
The universal cover, $\tilde{N}$, of $N$,
is a complete noncompact nilpotent Lie group.
Note $\pi_1(M)=\pi_1(N)$.} \cite{Wei}

{\bf Examples of Wilking (2000)}: 
{\em For any finitely generated
almost nilpotent group, $H$, one
can construct $M^n$ with $\pi_1(M^n)=H$.
This construction is done using $\tilde{N}$ from
Wei's construction, taking its
k fold isometric product and crossing with SU(2)
before dividing by $H$ and taking a similar warped
product.} \cite{Wlk}

Note that noncontractible loops in all these examples
slide to infinity.  This led the second author
to define the following concept in \cite{Sor4}.

\begin{figure}[htbp] 
\begin{center}
\includegraphics[height=1.5in ]{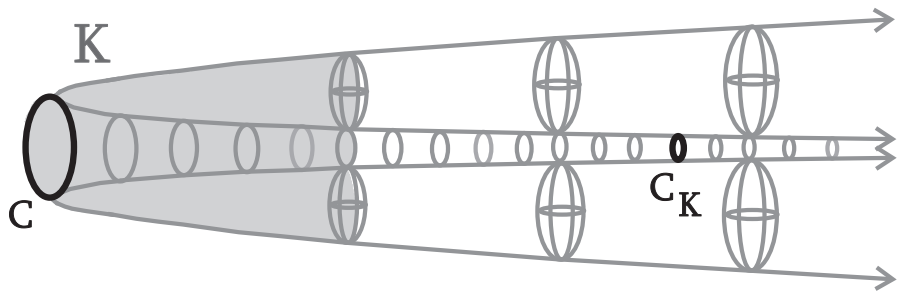}
\end{center}
\caption{Nabonnand's example
has the loops to infinity property as can be seen
taking the shaded compact set $K$.} \label{nab3} 
\end{figure}

\begin{figure}[htbp]
\begin{center}
\includegraphics[height=1in ]{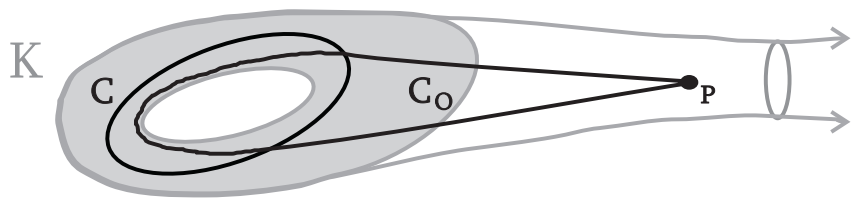}
\end{center}
\caption{
A complete punctured torus however,
does not have the loops to infinity
property because loops get caught on the
finite hole.
} \label{puncttorus}
\end{figure}

\begin{figure}[htbp]
\begin{center}
\includegraphics[height=.5in ]{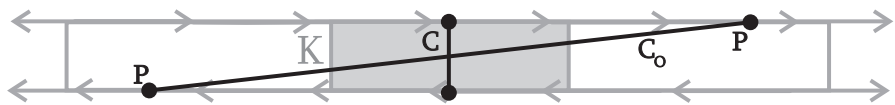}
\end{center}
\caption{
An infinite Moebius strip, defined as $\EE^2$
with points $(x,y)$ identified with $(-x, y+1)$,
also fails to have the loops to infinity property.  
} \label{moebius}
\end{figure}

\begin{figure}[htbp]
\begin{center}
\includegraphics[height=1in ]{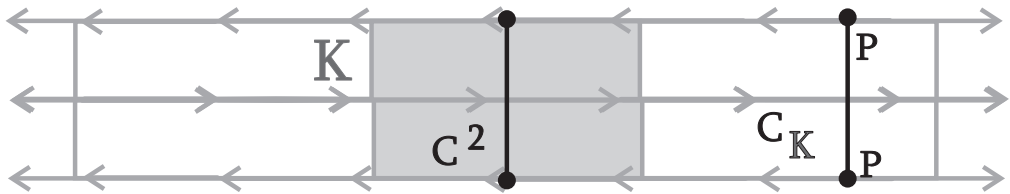}
\end{center}
\caption{However the infinite Moebius strip's double cover, the cylinder, satisfies
the loops to infinity property.
} \label{moebius3}
\end{figure}

A manifold has the {\bf loops to infinity} property if given any
noncontractible closed curve, $C$, 
and given any compact set $K$, 
then there exists a curve, $C_K$
contained in $M\backslash K$ which is
freely homotopic to $C$.  See Figures~\ref{nab3},
\ref{puncttorus}, \ref{moebius} and~\ref{moebius3}.

{\bf Sormani Loops to Infinity Theorem:}
{\em Either $M^n$ satisfies the 
loops to infinity property
or $M^n$ has a double cover 
which splits isometrically.} \cite{Sor4}

Clearly $M^n$ with
$Ricci>0$ must then have the loops to
infinity property.  
Two examples of $M^n$ with
split double coves are the
infinite Moebius strip depicted
in Figure~\ref{moebius}
and the oriented normal 
bundle over $\RR\textrm{P}^2$.

 {\bf Proof Outline:}
If there is a curve $C$ and a compact set $K$
such that any loop freely homotopic to $C$ passes
through $K$, then one takes a ray $h(t)$ and
$t_i \to \infty$ and the shortest loop $C_i$
freely homotopic to $C$ based at $h(t_i)$ must
pass through $K$.  
The lifts of $C_i$ to $\tilde{C_i}$ in
the universal cover, $\tilde{M}$ can be shown to
be minimal geodesics with length, $L(\tilde{C}_i)\to \infty$.
Transforming them by $g_i$ to a common compact lift of
$K$, one sees that 
a subsequence of the $g_i\tilde{C}_i$ converge to a line.
The Splitting Theorem 
then implies that $\tilde{M}$
splits.  It takes further work to split 
a double cover isometrically \cite{Sor4}.

{\bf Open Problem: The Fundamental Group at Infinity}
{\em in the sense of \cite{GeoMih} might 
contain $\pi_1(M)$.   The definition of the fundamental
group at infinity requires a proper homotopy yet in \cite{Sor4} 
the homotopy running from $C$ to $C_K$ as $K$ grows 
was not controlled.  }

It is a consequence of the Loops to Infinity property
that if $D$ is a compact domain in $M^n$
with one simply connected boundary, then $D$
is simply connected.  Earlier restrictions
on the fundamental group of such $D \subset M^n$ were
found by {\bf Schoen-Yau} in \cite{SchYau2} using {\em harmonic maps}.

Using the relative homology of arbitrary compact domains
in $M^n$ arising from the loops to infinity property and techniques
from algebraic topology
the authors proved the following:

{\bf Shen-Sormani (2000):} 
{\em Either 
$M^n$ is a flat normal bundle over a compact totally 
geodesic submanifold, 
or
$M^n$ has a trivial codimension one homology, $H_{n-1}(M^n,Z)$,
and $H_{n-2}(M,Z)$ is torsion free.} \cite{ShnSor}

This complete classification of $H_{n-1}(M,Z)$
extends {\bf S.T. Yau}'s 1976 proof using harmonic
forms that $M^n$ with $Ricci >0$ have trivial $H_{n-1}(M,R)$
\cite{Yau2}.  The first author had results in
this direction using Morse Theory in \cite{Shn}
and {\bf Itokawa-Kobayashi} had partially classified 
$H_{n-1}(M,Z)$ using minimizing currents \cite{ItKo}.

The control on $H_{n-2}(M,Z)$ is much weaker that
the control on the codimension one homology.  In fact,
examples have been constructed where $H_{n-2}(M,Z)$  
is infinite dimensional.  See
Section~\ref{examples}.

\sect{Volume and the Fundamental Group} \label{volume}

In this section we describe the properties of the fundamental
group of our open manifold, $M^n$, with nonnegative Ricci curvature.  
We focus on Milnor's unsolved conjecture:

{\bf Milnor Conjecture 1968}
{\em The fundamental group, $\pi_1(M^n)$, is finitely generated. 
That is, there are only finitely many one dimensional holes.} \cite{Mil}
Here we will survey partial results and obstructions
towards finding a possible counter example.  

The most useful tool for studying
the fundamental group of $M^n$ is its universal cover, $\tilde{M}^n$,
which is also open and has $Ricci \ge 0$.  Recall
the fundamental group acts on the universal cover by isometries called
``deck transforms'' and that $M^n= \tilde{M}^n / \pi_1(M)$.
If $i_0>0$
is the injectivity radius of $M$ then
\be \label{disjoint1}
d_{\tilde{M}}(gp, hp) \ge 2i_0 \qquad \forall g,h \in \pi_1(M), 
\textrm{ and } \forall p \in \tilde{M}.
\ee
The size and/or growth of the fundamental group can thus
be studied by counting disjoint balls of radius $i_0$
in the universal cover.  

For this reason, the Bishop-Gromov Volume Comparison Theorem plays a crucial
role in the study of $\pi_1(M^n)$ applied to both $M^n$ and $\tilde{M}^n$.

{\bf Bishop-Gromov Volume Comparison Theorem}

\noindent
{\em If $M^n$ has $Ricci \ge 0$, $p\in M^n$ and $0<r<R$ then
\be
\frac{Vol(B_p(r))}{Vol(B_p(R))} \ge \frac{\omega_n r^n}{\omega_n R^n}
\ee
where $\sigma_n=Vol(B_0(1)\subset \EE^n$.}

This theorem was proven by {\bf Gromov (1981)}  \cite{Gr2} using 
an estimate of {\bf Bishop (1963)} \cite{Bi}.
One can intuitively think of volumes as
capturing the fact that while some geodesics may bend apart
as they emanate from $p$ others must compensate by bending
together, thus the region they sweep out bends inward causing
inner balls to have larger volumes than expected compared to outer balls.
As a consequence $M^n$ has at most Euclidean volume growth:
\be
\limsup_{r\to\infty} \frac{Vol(B_p(r))}{r^n} \le \omega_n<\infty.
\ee
In fact $Vol_p(R) \le \omega_n R^n$ for all $R>0$.  

\begin{figure}[htbp]
\begin{center}
\includegraphics[height=3in ]{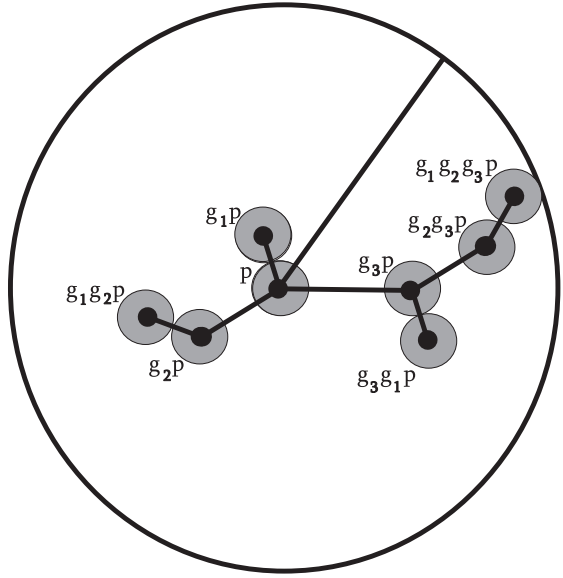}
\end{center}
\caption{Milnor's Estimate for $N(3)$ with $H\,=\,<g_1,g_2,g_3>$} 
\label{milnor}
\end{figure}

Milnor noticed that if one lists a finite collection of generators
\be
\{g_1, g_2, ... g_k\} \subset \pi_1(M) \textrm{ and let } H\,=\,<g_1,...g_k>
\ee
then $H$ cannot be too large or grow too
quickly.  If it did too many disjoint balls $B_{gp}(i_0)$ would fit
in $\tilde{M}$.  Particularly, {\bf Milnor} proved that 
\be
N(k)=\textrm{ the number of words of length } k\textrm{ in } H
\ee
grows at most polynomially in $k$ of order $n=dim(M^n)$ \cite{Mil}.
See Figure~\ref{milnor}.  

Each ball in Figure~\ref{milnor} 
has the same volume, due to the isometries, and each is disjoint
from the other by (\ref{disjoint1}).   They all fit in a large
ball of radius 
\be
R= k m \textrm{ where }m= \max_{g=g_1,g_2,g_3} d_{\tilde{M}}(p,gp).
\ee
because 
\be
d_{\tilde{M}}(g_1g_2p,p)\le
d_{\tilde{M}}(g_1g_2p,g_2p)+
d_{\tilde{M}}(g_2p,p)\le 2 m.
\ee
Here $m$ is finite because there is a finite list of generators.
To estimate $N(k)$, one sums over all words, $h$, of length $\le k$:
\be
N(k) {Vol(B_p(i_0/2))}= \sum_{h} Vol(B_{hp}(i_0/2)) \le \omega_n R^n
\le \omega_n (k \cdot m)^n.
\ee 
In Figure~\ref{milnor}, $H$ is abelian to make it easier to draw.
In general, $M^n$ need not have an abelian fundamental group.  

{\bf Gromov (1981)} proved that any finitely generated group of
polynomial growth is almost nilpotent \cite{Gr1}.  Combined with Wilking's
example described in Section~\ref{examples}, this would completely
classify the fundamental groups of $M^n$ if Milnor's Conjecture holds.

We now turn to partial solutions of the Milnor
Conjecture.  The most general is Wilking's result proven using
algebraic methods based on the qualities of $\pi_1(M)$
described above.

{\bf Wilking's Milnor Conjecture Reduction:} {\em If there exists
a counter example to the Milnor Conjecture then it
has a covering space with an abelian fundamental group
which is also infinitely generated.} \cite{Wlk}

All the other partial solutions involve additional
conditions on volume or diameter growth.

{\bf Maximal Volume Growth}

Recall that
the Bishop-Gromov Volume Comparison Theorem implies
$$
\limsup_{r\to\infty} \frac{Vol(B_p(r))}{r^n} \le \omega_n<\infty.
$$

{\bf Peter Li's Euclidean Volume Growth Theorem}: 
{\em If $M^n$ has Euclidean volume growth:
\be \label{maxvol}
\liminf_{r\to\infty} \frac{Vol(B_p(r))}{r^n} >0
\ee
then the fundamental group is finite.} \cite{Li}

The proof uses the heat kernel 
on the universal cover, $\tilde{M}$.  
It is interesting to note
that in dimension $3$, $H_1(M^3,Z)$ is torsion free by
\cite{ShnSor}, thus the only three dimesional manifold with
Euclidean volume growth has $H_1(M,Z)=0$.  In fact Zhu
proved $M^3$ satisfying \ref{maxvol} is contractible \cite{Zhu1}.

{\bf Anderson's Volume Growth Theorem:} 
{\em If $M^n$ has $b_1(M)\ge k$ and
$$
\limsup_{r\to\infty} \frac{Vol(B_p(r))}{r^{n-k}} >0
$$
then $ \pi_1(M)$ is finitely generated.} \cite{And1}

This implies Li's result when $k=0$ and is
proven using a clever volume comparison arguement
relating large balls in $M$ to their lifts
in $\tilde{M}$ restricted to fundamental domains.
Anderson  also obtains estimates on $b_1(M)=dim(H_1(M,Z)$ assuming
additional sectional curvature bounds \cite{And1}.

{\bf Perelman Contractibility Theorem}: 
{\em $M^n$ is contractible if it almost maximal volume growth:
$$
\liminf_{r\to\infty} \frac{Vol(B_p(r))}{r^n} >\omega'_n
$$
where $\omega'_n$ is sufficiently close to $\omega_n=vol(B_0(1)\subset \EE^n)$.
} \cite{Per1}

In Section~\ref{examples}, we present Menguy's four dimensional example
demonstrating $\omega'_n$ must be close to $\omega_n$
for this result to hold.  
Zhu proved that three dimensional $M^3$ satisfying only (\ref{maxvol})
are contractible \cite{Zhu1}.

Note that for a presumably larger constant, $\omega'_n$,
also depending only on dimension,  {\bf Cheeger-Colding}
have proven that $M^n$ is diffeomorphic to $\RR^n$.  Their
proof uses almost rigidity techniques and so one cannot estimate
the actual value of their constant \cite{ChCo} Thm A.1.12.  We cannot
describe their proof in the space allowed here but do
describe Perelman's:

{\bf Proof Outline:}
Perelman proves that $f:S^k \to M$ can be continuously
extended to $f:D^{k+1} \to M$ using induction on $k$.
His key estimate depending on volume is proven using 
Bishop-Gromov's proof of the volume comparison theorem.
It says that for any $c_2>c_1>0$ and $\epsilon>0$
there exists $\delta>0$ such that if 
\be \label{pervol}
Vol(B_p(c_2R))>(1-\delta)(c_2R)^n
\ee
then for every $a \in B_p(c_1R)$ there exists
$b \in M\setminus B_p(c_1R)$ such that the geodesic
from $p$ to $b$ passes through $B_a(\epsilon R)$.
If this were not true there would be too many
cut points and the volume of the larger ball would
not be almost maximal contradicting (\ref{pervol}).
This allows Perelman to proceed with a filling in 
procedure for each subsequent cell $k$.  See \cite{Per1}
for illustrations and details.  By analyzing the
process carefully one can estimate the value of $\omega'_n$.

The key intuition here is that homology causes cut points
and cut points use up space.

{\bf Minimal Volume Growth}

S-T Yau \cite{Yau1} used results on harmonic functions to prove 
$M^n$ has at least linear volume growth:
$$
\liminf_{r\to\infty} \frac{Vol(B_p(r))}{r}=C_{M^n} > 0.
$$
This can be reproven by applying the Bishop-Gromov volume comparison
to balls around points along a ray.

The second author proved that if $M^n$ has at most linear
volume growth
$$
\limsup_{r\to\infty} \frac{Vol(B_p(r))}{r} <\infty
$$
then it has sublinear diameter growth
$$
\limsup_{r\to\infty} \frac{diam(\partial B_p(r))}{r} =0,
$$ 
using Cheeger-Colding almost rigidity techniques.
This diameter is extrinsic so $diam(\partial B_p(r))\le 2r$. \cite{Sor1}\cite{Sor2}

{\bf Sormani Small Linear Diameter Growth Theorem}: 
{\em If $M^n$ has 
small linear diameter growth:
$$
\limsup diam(\partial B_p(r))/r < S_n
$$
then $M^n$ has a finitely generated 
fundamental group. } \cite{Sor3}

{\bf Corollary:} {\em The Milnor Conjecture 
is proven for manifolds, $M^n$,
with minimal volume growth.}

The constant $S_n$ was given explicitly 
in \cite{Sor3} and then was improved in 2003 by
S Xu, Z Wang and F Yang \cite{XuWaYa}.
Recently, W. Wylie observed that $\pi_1(M^n)$
is finitely presented when $M^n$ has
small linear diameter growth \cite{Wy}.

{\bf Proof Outline:} Since $M^n$ is complete,
one can construct special {\em halfway generators}, 
$g_i \in \pi_1(M^m, p),$
with loops $c_i:[0,L_i]\to M$ such that
\be \label{halfeq}
d(c_i(0), c_i(L_i/2))=L_i/2.
\ee

\begin{figure}[htbp]
\begin{center}
\includegraphics[height=1.2in ]{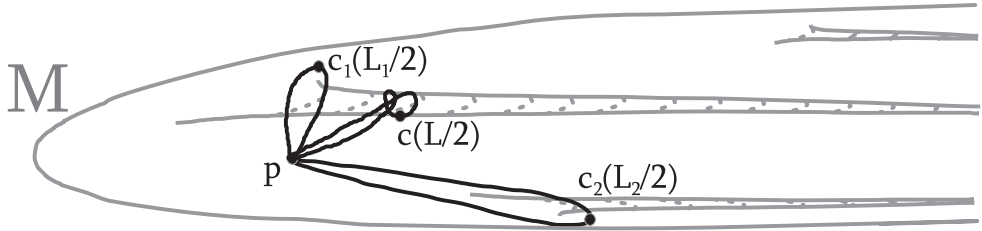}
\end{center}
\caption{Halfway generators: $c_i$ satisfy (\ref{halfeq}), $c$ does not.} 
\label{Sormani1}
\end{figure}

Using $Ricci\ge 0$ Sormani proves the
halfway generators' loops, $c_i$, 
satisfy a {\em Uniform Cut} property on
a ball $B_i=B_{c_i(L_i/2)}(S_nL_i)$:
{\em all geodesics from $p$ entering $B$ are cut}.

\begin{figure}[htbp]
\begin{center}
\includegraphics[height=1in ]{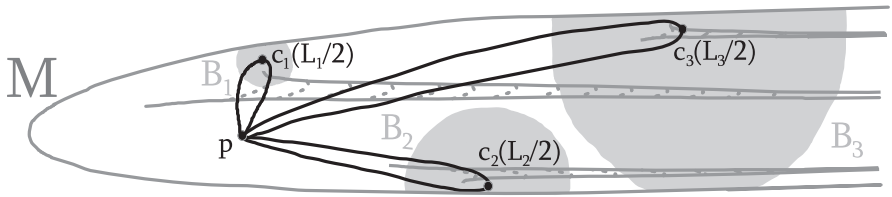}
\end{center}
\caption{All geodesics from $p$ entering $B_i$ are cut.
Thus there is no room for a ray to pass through $B_3$.} \label{Sormani2}
\end{figure}

This contradicts the existence of a
ray in a complete noncompact manifold
when the diameter growth is too small.

The {\em Uniform Cut} property on the ball 
$$
B=B_{c_i(L_i/2)}(S_nL_i)
$$
is proven by applying the
Abresch-Gromoll Excess Theorem 
to the lift of the $c_i$ to $\tilde{M}$. See \cite{Sor3}
for details.

In the next section we will see that we cannot hope
to get contractibility for manifolds with minimal
volume growth as Menguy has constructed an example
with bounded diameter growth and infinite topological
type.  

Intuitively the space used up by cut points
resulting from 1 dimensional holes \cite{Sor3}
is significantly larger
than the space used up by cut points resulting from
higher dimensional holes \cite{Per1}.

\sect{Examples} \label{examples}

In this section we describe a host of important examples
of open manifolds with nonnegative Ricci curvature
supplementing those of Nabonnand, Wei and Wilking described
in Section~\ref{geods}.  
The first manifold with positive Ricci curvature and infinite
topological type was constructed by Sha and Yang.  It was seven
dimensional with infinite dimensional $H_4(M,Z)$.  \cite{ShaYng1}
This result was then generalized in 1991 as follows:

{\bf Sha-Yang Examples}:
{\em For all integers $p\ge 2$ and $q\ge 1$
there exists $M^{p+q}$ with nonnegative Ricci
curvature which is created from
$S^{p-1} \times R^{q+1}$ by cutting off
infinitely many $S^{p-1}\times D^{q+1}_i$ and 
gluing in $S^p \times S^q_i$. 

In particular, there exists a 4 dimensional
manifold $M^4$ with infinite second Betti number and $Ricci >0$.
\cite{ShaYng2}}

\begin{figure}[htbp]
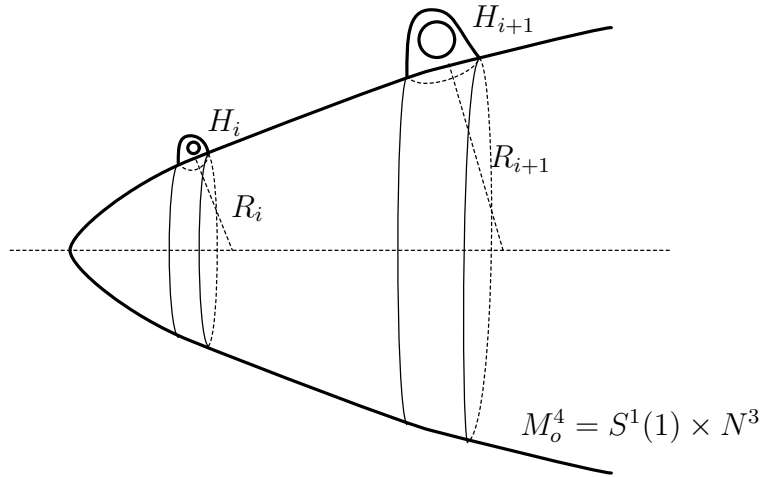

\begin{center}
\begin{texdraw}
\drawdim cm
\linewd  0.04
\relsegscale 0.8

\move(0 0)
\clvec(0 0.3)(1 1.1)(2 1.5)
\clvec(3 1.9)(5.6 2.9)(6 3)
\clvec(6.4 3.1)(8.8 3.7)(9 3.7) 

\move(1.8 1.45)
\clvec(1.8 1.7)(1.8 1.9)(2 1.9)
\clvec(2.2 1.9)(2.3 1.7)(2.3 1.6)

\move(2.05 1.7)\lcir r:0.1

\move(0 0)
\clvec(0 -0.3)(1 -1.1)(2 -1.5)
\clvec(3 -1.9)(5.6 -2.9)(6 -3)
\clvec(6.4 -3.1)(8.8 -3.7)(9 -3.7)

\move(5.6 2.9) 
\clvec(5.6 3.5)(5.6 4)(6 4)
\clvec(6.4 4)(6.4 3.7)(6.8 3.2)

\move(6.1 3.5)
\lcir r:0.3

\linewd 0.02

\lpatt(0.05 0.05)
\move(2.1 1.5)\lvec(2.7 0)

\move(2.7 0.7)
\textref h:L v:C \htext{$R_i$}

\move(-1 0)\lvec(10 0)
\move(7.2 0)\lvec(6.3 3.1)

\move(1.8 1.45)\clvec(2 1.2)(2.3 1.4)(2.3 1.6)
\move(2.3 1.6)
\clvec(2.5 1.4)(2.5 -1.4)(2.3 -1.6)

\move(5.6 2.9)\clvec(5.8 2.6)(6.7 2.9)(6.8 3.2)
\move(6.8 3.2)
\clvec(7.1 2.9)(7.1 -2.9)(6.6 -3.2)

\lpatt()
\linewd 0.02
\move(1.8 1.45)
\clvec(1.6 1.25)(1.6 -1.25)(1.8 -1.45)
\move(2.3 1.6)
\clvec(2.1 1.4)(2.1 -1.4)(2.3 -1.6)

\move(5.6 2.9)
\clvec(5.4 2.6)(5.4 -2.6)(5.6 -2.9)
\move(6.8 3.2)
\clvec(6.5 2.9)(6.5 -2.9)(6.6 -3.2)

\move(7.0 1.5)
\textref h:L v:C \htext{$R_{i+1}$}
\move(7.5 -3.2)

\textref h:L v:B \htext{$M^4_o=S^1(1) \times N^3$} 

\move(2.3 1.9)
\textref h:L v:B \htext{$H_i$}
\move(6.7 3.6)
\textref h:L v:B \htext{$H_{i+1}$}

\end{texdraw}
\end{center}
\caption{The 4D Sha-Yang Example with $p=2$ and  $q=2$} 
\label{shayangfig}
\end{figure}

Sha and Yang focus on the compact case in this paper only briefly
sketching the construction of the open manifold we depict
in Figure~\ref{shayangfig}.  More detail on the open
manifold is available in \cite{ShnWei} where $M^4$  is shown to have
$diam(\partial B_p(r) \le C r^{3/4}$ and $vol(B_p(r)) \le r^{5/2}$. 

{\bf Construction:}
The construction when the dimension $n=4$
begins with  $M^4_o = S^1 (1) \times N^3$ with 
the isometric product metric, where $N^3\approx R^3$
When one stays away from  
the tip of $M^4_o$, the size of $S^1(1)$ is relatively small, compared 
to the cross-sections of $N^3$. Thus 
in the Figure~\ref{shayangfig} for $M^4_o$, we make $M^4_o$ look like $N^3$.
Note $M_0^4$ has $Sect \ge 0$. 

$N^3$ is a surface of 
revolution in $R^4$ which looks something like a paraboloid
but has a sequence of annual regions, $A_i$, with constant sectional
curvature $1/R_i^2$ and width $2r_i=2 \alpha_iR_i$.  This can be achieved
by taking a hemisphere in $S^3$ and slicing it into infinitely many
annuli of width $2 \alpha_i$ such that 
$\sum_{i=1}^{\infty}2\alpha_i < \pi/2$, then spreading these
annuli apart from each other and rescaling each up by some
radius $R_i$ to create an open manifold $N^3$.

To create the interesting topology, Sha-Yang edit $M_0^4$.
From each $A_i$ which is isometric to an annular region in $S^3$
crossed with $S^1$, they remove a metric ball 
$B_i^3 := B_i(p_i, r_i)$ from $N^3$.   This creates a manifold 
\be
\hat{M}^4 := S^1(1)\times (N^3\setminus \coprod_{i=1}^{\infty} B^3_i )
\ee
with infinitely many boundaries. Each boundary is an isometric
product of an $S^1 \times \partial B_i^3 \subset S^1 \times_{R_i} S^3$
where $R_i$ is just the constant radius.

Topologically, they can glue a handle  $H:=D^2(1)\times S^2(1)$   
to $\hat{M}^4$ along these  boundaries. The resulting manifold is a 
manifold $M^4$ with $H_2(M^4, Z)$ infinitely generated: 

\[ M^4 := S^1(1)\times ( N^3\setminus\coprod_{i=1}^{\infty} B^3_i) \bigcup_{Id}  \coprod_{i=1}^{\infty} H_i\]
where each $H_i$ is diffeomorphic to $H$.  See Figure~\ref{shayangfig}.

In order to understand why $M^4$ has positive Ricci curvature, we
now describe the doubly warped product on the handles 
\be
H_i = (0, r_i) \times_{h_i} S^1(1) \times_{f_i} S^2(1). 
\ee  
To obtain the handle topology we desire while closing smoothly
at $r=0$ we set
$h_i(0)=0$, $f_i'(0)=0$ and $f_i(0)=h_i'(0)=1$.
To smoothly attach this into $\partial \hat{M} \cap A_i$
we need to attach the $S^1$ direction
straight (so $h_i(r_i)=1$, $h_i'(r_i)=0$)
and the other directions curved like a sphere of radius $R_i$
(so $f_i(r) =R_i  \sin ( r/R_i)$ for $r$ nearby $r_i$).
Sha-Yang then carefully choose $f$ and $g$ satisfying
(\ref{eqnwarp}) and the other equations guaranteeing positive
Ricci curvature.   \cite{ShaYng2}

In Sha-Yang's surgery, a removed part, $S^{p-1}\times D^{q+1}$, can be
viewed as a trivialization of the normal bundle of $S^{p-1}$ in
$M^{p+q}:=S^{p-1}\times N^{q+1}$ with the metric product metric, where
$N^{q+1}$ can be a round sphere or a surface of revolution. David {\bf Wraith}
studies the  surgery problem on a manifold  $M^{p+q}$ with positive
Ricci curvature and surgeries of codimension three. His technique 
is similar to Sha-Yang's.  Wraith  needs
the same local form for the metric on the ambient manifold $M^{p+q}$ in
order to complete the Ricci positive surgery. The essential difference
is that he handles the surgery with a non-standard trivialization which is
not determined by the metric. To do so, he assumes that $p\geq q+1\geq
3$ in order to use a smooth map $T:  S^{p-1} \to SO(q+1)$ to make a
twisting for 
\be
\tilde{T}: (x, r, y)\in S^{p-1} \times D^{q-1}\to  
(x, r, T(x)y)\in S^p \times D^{q+1}
\ee
 when gluing in 
$D^p\times S^q$. Apart from the restriction on
dimensions, Wraith's technique can be used in all situations where
Sha-Yang's technique can be used.  It can also be applied to exotic spheres.
Of interest here is that he can construct 
complete open manifolds with positive Ricci curvature and infinite
topological type with are similar to the Sha-Yang examples
but not diffeomorphic to them. \cite{Wra}

{\bf Anderson-Kronheimer-LeBrun Examples (1989):}
{\em $M^4$ which are Ricci flat and Kahler 
with infinite dimensional $H_2(M,Z)$
based on physics of Gibbons-Hawking}.\cite{AndKrLb}

{\bf Construction:}
Anderson-Kronheimer-Lebrun's example is constructed from the 
Gibbons-Hawking Ansatz by using infinitely many, sparsely 
distributed centers. It requires some more expertise
than the other examples to understand.

Take a sequence of points $p_j = (j^2, 0,  0)$  in $R^3$. 
There is a unique principal $S^1$-bundle 
$\pi_o: M_o \to R^3_o:= R^3 \setminus \{ p_j \}$ 
such that the Chern class is $-1$ when restricted to a  
sphere $S^2(p_j, r_j)\subset R^3$ for small 
$r_j < \min_{k\not=j} \|p_j - p_k \|$. 
An important fact is that $\pi_0^{-1} ( B(p_j,r_j )) $ 
is diffeomorphoc to $ \hat{B}^4_j:= B^4_j \setminus\{ 0\}$,  
where $B^4_j  $ is a copy of a ball in $R^4$, such that the 
action of $S^1\subset \mathbb C$ on 
$\hat{B}^4_j  \subset  R^4 \approx \mathbb C^2$ 
is given by scalar multiplication. 

\begin{figure}
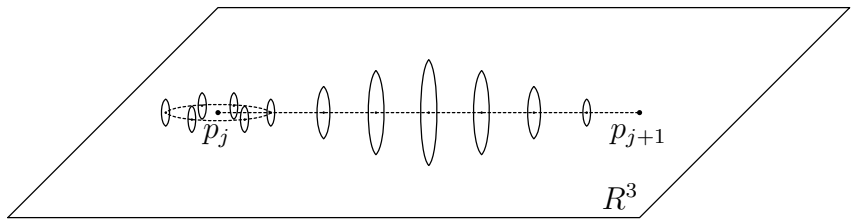

\begin{center}
\begin{texdraw}
\drawdim cm
\linewd  0.02
\relsegscale 0.7

\move(-8 -2)\lvec(4 -2)\lvec(8 2)\lvec(-4 2)\lvec(-8 -2)
\move(-4 0)\fcir f:0 r:0.05
\move(4 0)\fcir f:0 r:0.05

\move(0 0)\fcir f:0 r:0.02
\move(0 -1)\clvec(-0.2 -0.8)(-0.2 0.8)(0 1)
\move(0 -1)\clvec(0.2 -0.8)(0.2 0.8)(0 1) 

\move(-1 0)\fcir f:0 r:0.02
\move(-1 -0.8)\clvec(-1.2 -0.6)(-1.2 0.6)(-1 0.8)
\move(-1 -0.8)\clvec(-0.8 -0.6)(-0.8 0.6)(-1 0.8)

\move(1 0)\fcir f:0 r:0.02
\move(1 -0.8)\clvec(1.2 -0.6)(1.2 0.6)(1 0.8)
\move(1 -0.8)\clvec(0.8 -0.6)(0.8 0.6)(1 0.8) 

\move(-2 0)\fcir f:0 r:0.02
\move(-2 -0.5)\clvec(-2.16 -0.35)(-2.16 0.35)(-2 0.5)
\move(-2 -0.5)\clvec(-1.84 -0.35)(-1.84 0.35)(-2 0.5) 

\move(2 0)\fcir f:0 r:0.02
\move(2 -0.5)\clvec(2.16 -0.35)(2.16 0.35)(2 0.5)
\move(2 -0.5)\clvec(1.84 -0.35)(1.84 0.35)(2 0.5) 

\move(-3 0)\fcir f:0 r:0.02
\move(-3 -0.25)\clvec(-3.1 -0.2)(-3.1 0.2)(-3 0.25)
\move(-3 -0.25)\clvec(-2.9 -0.2)(-2.9 0.2)(-3 0.25)

\move(3 0)\fcir f:0 r:0.02
\move(3 -0.25)
\clvec(3.1 -0.2)(3.1 0.2)(3 0.25)
\move(3 -0.25)\clvec(2.9 -0.2)(2.9 0.2)(3 0.25)

\move(-5 0)\fcir f:0 r:0.02
\move(-5 -0.25)\clvec(-5.1 -0.2)(-5.1 0.2)(-5 0.25)
\move(-5 -0.25)\clvec(-4.9 -0.2)(-4.9 0.2)(-5 0.25)

\move(-3.5 -0.13)
\fcir f:0 r:0.02
\move(-3.5 -0.38)
\clvec(-3.6 -0.33)(-3.6 0.09)(-3.5 0.12)
\move(-3.5 -0.38)\clvec(-3.4 -0.33)(-3.4 0.09)(-3.5 0.12)

\move(-4.5 -0.13)
\fcir f:0 r:0.02
\move(-4.5 -0.38)
\clvec(-4.6 -0.33)(-4.6 0.09)(-4.5 0.12)
\move(-4.5 -0.38)\clvec(-4.4 -0.33)(-4.4 0.09)(-4.5 0.12)

\move(-3.7 0.13)
\fcir f:0 r:0.02
\move(-3.7 0.38)
\clvec(-3.8 0.33)(-3.8 -0.09)(-3.7 -0.12)
\move(-3.7 0.38)\clvec(-3.6 0.33)(-3.6 -0.09)(-3.7 -0.12)

\move(-4.3 0.13)
\fcir f:0 r:0.02
\move(-4.3 0.38)
\clvec(-4.4 0.33)(-4.4 -0.09)(-4.3 -0.12)
\move(-4.3 0.38)\clvec(-4.2 0.33)(-4.2 -0.09)(-4.3 -0.12)

\lpatt(0.05 0.05)
\move(-4 0)\lvec(4 0)

\move(-3 0)\clvec(-3.2 0.2)(-4.8 0.2)(-5 0)
\move(-3 0)\clvec(-3.2 -0.2)(-4.8 -0.2)(-5 0)

\move(-4 -0.2)
\textref h:C v:T \htext{$p_j$}
\move(4 -0.2)
\textref h:C v:T \htext{$p_{j+1}$}

\move(4 -1.8)
\textref h:R v:B \htext{$R^3$}
\end{texdraw}
\end{center}
\caption{Anderson-Kronheimer-LeBrun Example} \label{akl}
\end{figure}

Then define 
\[ M^4 := M_o\bigcup_{Id} \coprod_{j=1}^{\infty}B^4_j.\]
This $M^4$ is a smooth manifold with 
$H_2 (M, Z) = \oplus_{j=1}^{\infty} Z $. 
$\pi_o$ can be extended to a map 
$\pi: M^4\to R^3$ such that $\pi^{-1}(p_j)$ is a point. 
Recall that $\pi^{-1}(q)$ is $S^1$ when $q \notin \{p_j\}$.
To illustrate $M^4$, we view $R^3$ as a plane in Figure~\ref{akl}. 

We now describe the construction of the metric
$g_0$ depicted in Figure~\ref{akl}.
The Chern class of $\pi_o : M_o \to R^3_o$ 
is represented by the closed 2-form 
$\frac{1}{2\pi}* df $. Let $\omega\in \Omega^1 (M_o)$ 
be a connection $1$-form for  $\pi_o: M_o\to R^3_o$ such that 
\[ \pi_o^*(* df ) = d \omega.\]
$\omega$ is unique up to a gauge transformation 
since $R^3_o$ is simply connected. The canonical metric 
on $M^4_o$ is defined by 
\[ g_o := \omega \otimes \omega + \pi_o^* d s_3^2,\]
where $ds^2_3$ denote the Euclidean metric on $R^3$. 
It is singular at the points $\pi^{-1}(p_j)$.

Next, Anderson-Kronheimer-LeBrun warp the metric $g_o$ as follows,
\[ g := f^{-1} \omega \otimes \omega + f \pi_o^* ds^2_3,\]
where 
$f: R^3_o \to R$ is defined by 
\[  f (x):= \frac{1}{2} \sum_{j=1}^{\infty} \frac{ 1}{ \| x - p_j \|}.\]
Anderson-Kronheimer and LeBrun then verify that 
\[ \Delta f = d * d f = 0.\]
This fact is true only in three dimensions!  Thus the metric $g$ 
continues smoothly across the isolated points $\pi^{-1}(p_j)$, so 
$g$ is a complete Ricci-flat metric.  \cite{AndKrLb}

{\bf Menguy's Euclidean Volume Growth Examples (2000):} 
{\em $M^4$ with $Ricci >0$, volume growth like $r^4$,
and infinite second homology.  } \cite{Mng1}

\noindent
The homology is created by cutting
and pasting in special convex manifolds with boundary constructed
by Perelman for a compact example in 1997.  We begin
by describing Perelman's construction \cite{Per2}.

\begin{figure}[htbp]
\begin{center}
\includegraphics[height=2in ]{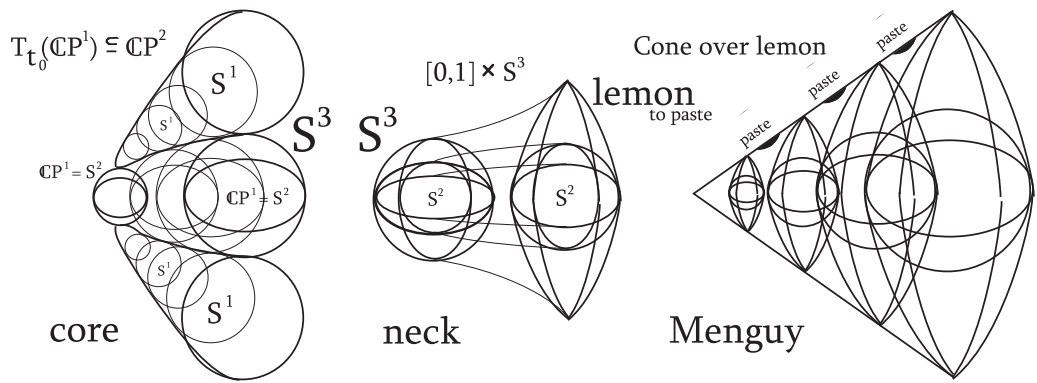}
\end{center}
\caption{Gluing the core to the neck, then repeatedly into Menguy's $M_0$.} 
\label{PerMen}
\end{figure}

{\bf Perelman's Building Blocks} have a core and a neck.  
The core introduces
the topology.  It views $S^3$ as an $S^1$ bundle over $S^2$, and
warps the $S^1$ down to $0$, while keeping the base $S^2$ positive,
so that the core has a noncontractible two sphere and a round
convex boundary $=S^3$.  Perelman justifies the smooth closing up of
the $S^1$ direction, by relating the warping direction to the
distance function from $CP^1=S^2$ in $CP^2$ giving a kind of cylindrical 
coordinates expression for the tubular neighborhood of the $CP^1$.
The boundary of this tubular neighborhood is an $S^1$ bundle over
$S^2$ which is diffeomorphic to the standard $S^3$.  It is convex
and all the geodesics are curving together as they approach the
boundary so it is difficult to glue a small copy of it into
a manifold with nonnegative Ricci curvature: all the geodesics
would be forced apart.    See Figure~\ref{PerMen}

In fact Colding's stability theorem \cite{Co}
states that it is impossible to glue in tiny pieces of topology
into a manifold with nonnegative Ricci curvature.  
Thus Perelman creates a neck which bends some geodesics outward
and others continue inward, so that he can glue the building block
into the singular edge of manifold.  The geodesics which bend
together fold over the edge and those that bend apart turn outward
along the edge and the whole manifold is then smoothed into
a smooth manifold with $Ricci >0$.

The neck is a doubly warped $S^3\times [0,1]$.  $S^3 \times\{0\}$ 
is a small round sphere that fits the core.
Then $S^3\times\{r\}$ grows towards a convex
boundary at $r=1$.  $S^3 \times \{1\}$ has a metric that looks
like a lemon.  It is a rotationally symmetric $S^3$ with 
the distance between the poles $=\pi R$ and a waist $=2\pi r$
where $r^{3-1}< R^3 <1$, so the sectional curvature is $>1$
and the normal curvatures are all $>1$.  

The core and the neck together form Perelman's building blocks which he
glued into a compact manifold close to a singular manifold formed by
taking a double spherical suspension over a round two sphere.  
When holes
are cut out of the singular edge they look like singular lemons, so
that when the singular manifold is smoothed slightly along the singular edge,
the holes can be made to precisely fit the lemon shaped boundary
of the building block.  To glue a tiny copy of the building block, the edge
must be sharper, reflecting Colding's Stability Theorem \cite{Co}.

{\bf Menguy's Construction:} 
Since Menguy's example is open, he needs
to insure that it is asymptotically singular at infinity in order
to successfully complete the editting process.  In fact Menguy starts
with a metric cone, $M_0$, over a spherical suspension of a small ball,
\be
dr^2 + (cr)^3 ( ds^2 + sin^2(s)R^2_0 d\sigma^2) \textrm{ where } c<1
\textrm{ and } R_0<1 
\ee
which has volume growth like $Cr^4$ and
two singular rays emanating from the pole at $r=0$.  He
smooths along these rays so that they have sharper and sharper
corners as $r\to \infty$.  Menguy is able to cut out a sequence of
smaller and smaller lemons and glue in the Perelman building blocks
so that his manifold has infinite second Betti number and 
volume growth  $Cr^4$.  
See Figure~\ref{PerMen}.  
Menguy also claims one could similarly
construct $M^{2n}$ with all even Betti numbers infinite although
this construction is only briefly outlined. \cite{Mng1}

{\bf Menguy's Bounded Diameter Example (2000):} 
{\em On the other extreme Menguy constructs an $M^4$ with bounded diameter
growth.  This is achieved using the same Perelman building blocks editted
into a manifold which is shaped somewhat like a sword:
\be
dr^2 + f(r)^2 ( ds^2 + cos^2(s)R_0 d\sigma^2),
\ee
where $f(r)$ is bounded above.}  \cite{Mng2}

\vspace{.2cm}
\noindent{ \bf Vector Bundles}

According to Cheeger-Gromoll's soul theorem, every complete open 
manifold of $K \geq 0$ is diffeomorphic to a vector bundle over 
a closed manifold with $K \geq 0$ \cite{ChGl1}. 
Thus it is a natural problem whether or not a vector bundle 
over a closed manifold with $K \geq 0$  admits a complete 
metric with $K \geq 0$. Although the topology of complete 
open manifolds with $Ric \geq 0$ is much more complicated,  
one can also study an analogue of the problem for  vector bundles 
over a closed manifold with $Ric \geq 0$. 

The first notable result in this direction is due to 
{\bf Nash} and {\bf Berard-Bergery}
(\cite{Nsh}, \cite{Ber}). They prove that 
every vector bundle of rank $\geq 2$ over a compact manifold with 
$Ric >0$ admits a complete metric with $Ric >0$. Note that rank 
one vector bundles have lines so they cannot admit
a metric with  $Ric >0$ by the Cheeger-Gromoll 
splitting theorem.   In fact Berard-Bergery proves that 
if $M^n$ has $Ricci\ge 0$ then there is a metric on 
$\RR^p \times M^n$ with $Ricci >0$ for all $p \ge 3$  \cite{Ber}.   
{\bf Otsu} constructed manifolds
with $Ricci\ge 0$ and  Euclidean volume growth
diffeomorphic to $\RR^{n-k} \times S^k$
and $\RR^{n-k}\times \RR P^2$ where $k \ge 2$ \cite{Ot}.
According to {\bf Anderson} \cite{And1}, no  $R^2$-bundle over 
a torus admits a complete metric with $Ric >0$. 

In 2002, {\bf Belegradek-Wei} constructed metrics of positive Ricci curvature
on vector bundles over nilmanifolds.  For any complex line bundle, $L$, 
over a nilmanifold, there is a sufficiently large $k$ such that the 
the Whitney sums of $k$ copies of $L$ admit metrics of positive Ricci 
curvature.  In particular if
$B$ has $ Ricci \ge 0$ and $\dim T\ge 4$, then for any sufficiently 
large $k$, there are infinitely many rank $k$ vector bundles over 
$B\times T$ with topologically distinct total spaces which admit 
metrics of $Ricci>0$ but are not homeomorphic to manifolds of 
$\sec\ge 0$ by work of Belegradek-Kapovitch \cite{BelWei1}\cite{BelKap}.

In 2004, {\bf Belegradek-Wei} construct further examples as follows.
Let $B$ be a closed manifold with $Ric\geq 0$. 
If $E(\xi)$ is the total space of a vector bundle $\xi$ over $B$, 
then $E(\xi)\times \RR^p$ admits a complete Riemannian metric 
with $Ric >0$ for all large $p$.
In fact, $B$ need not have $Ricci \ge 0$ as long as $B$ is the 
total space of a smooth fiber bundle $ F \to B\to S$, where 
$F,S$ are  closed manifolds with $Ric \geq 0$ and the 
structure of the bundle lies in the isometry group of $F$. 
It can be shown that such a manifold $B$ admits a metric 
of almost nonnegative Ricci curvature, but has no
metric of $Ricci \ge 0$ yet $E(\chi)\times \RR^p$ will still admit 
a metric with positive Ricci curvature. 
The minimum value of $p$ depends on $B$ and its bundle structure. 
The constructions are warped products 
$\RR^+ \times_f E(\chi) \times_h S^{p-1}$ where $f(0)=1$ and $h(0)=0$.
For further details and many explicit examples see \cite{BelWei2}.

For additional information about open manifolds with sectional
curvature $K \ge 0$ see Greene's survey \cite{Gre}.

\sect{Three Manifolds} \label{3D}
In this section we review the properties of complete noncompact 
3 manifolds, $M^3$, with nonnegative Ricci curvature.  The most
substantial contribution to this topic is in Schoen-Yau's 1982
paper:

{\bf Three Manifold Theorem of Schoen-Yau (1982)}
{\em
If $M^3$ has $Ricci >0$ then $M^3$ is diffeomorphic to $\RR^3$.}\cite{SchYau1}

{\bf Open Problem on Three Manifolds with {\em Ricci} $\ge$ 0:}
{\em Classify the topology of $M^3$ 
and prove the Milnor Conjuecture in dimension three.}

Schoen and Yau never had the opportunity to investigate
this nontrivial problem further after their original
paper on manifolds with positive Ricci curvature.  Recently Schoen
has suggested that {\em Ricci flow} might be used to prove that
these manifolds are either diffeomorphic
to manifolds with split covering spaces or manifolds
with positive Ricci curvature in which case their
theory should apply.    However when {\bf Shi} tried 
using Ricci flow he needed an additional upper bound
on sectional curvature to ensure the uniqueness of
the flow \cite{Shi}.  Perhaps more recent methods on Ricci flow would
prove effective.  

Here we will describe what is known about this open problem in
some detail and provide some relevant information from three
manifold topology.  We begin by reviewing Schoen-Yau's paper.

{\bf Proof Outline:}
Schoen-Yau begin with a proof that $M^3$ {\em with only $Ricci \ge 0$
has $\pi_2(M)=0$ unless $M$ is isometrically covered by a product
of a real line with a compact surface $S^2$} which involves
the Cheeger-Colding splitting theorem.  
To study the fundamental
group, they then assume $\pi_2(M)=0$, that $M$ is orientable and the 
fundamental group is $\ZZ$ generated by an element represented
by a curve $\sigma$.  They use Poincare Duality to create a sequence of
compact orientable surfaces, $\Sigma_i$, which intersect $\sigma$,
such that $\partial \Sigma_i \subset \partial M_i$ and $\bigcup M_i=M$.
They, in fact, choose these $\Sigma_i$ to be minimal in their
homotopy class, thus allowing them to use minimal surface methods to
prove that a subsequence of the $\Sigma_i$ converge to a complete
{\em noncompact stable minimal surface} $\Sigma$.  
Fischer-Colbrie and Schoen had proven
earlier that such a $\Sigma$ {\em must be totally geodesic and $Ricci(N,N)=0$
where $N$ is any normal to $\Sigma$}. 

{\bf Adding the assumption that {\em Ricci} $>$ 0}, Scheon-Yau 
conclude that $M^3$ is contractible since such $\Sigma$ cannot exist.
They then proceed to prove $M^3$ is simply connected
at infinity as well, again using a minimal surface contradiction
arguement involving positive Ricci curvature.  Finally they prove
$M^3$ is irreducible using a complicated technical arguement that
Kleiner notes would be greatly simplified if one assumes that
Poincare Conjecture has been proven.  Applying a result of Stallings
they complete the proof of their theorem.

{\bf Assuming only {\em Ricci} $\ge$ 0}, one would have to carefully examine 
the possibilities that arise from the existence of such a minimal surface 
$\Sigma$.     {\bf Meeks-Simon-Yau} \cite{MSY} 
have done partial work in this direction
showing a compact three manifold with mean convex boundary 
that has $Ricci \ge 0$
is a solid handlebody.   Note {\bf Ananov-Burago-Zalgaller} \cite{AnBuZg} 
obtained the same result by studying the Morse theory of
the distance function to the boundary rather than $\Sigma$.  
{\bf Anderson-Rodriguez}
 \cite{AndRod} added the additional assumption that sectional
curvature is bounded above and proved that the existence of $\Sigma$
implies $M^3= \Sigma \times \RR^+$ when $\Sigma$ and $M$ are oriented.  
This agrees with Shi's results using Ricci flow \cite{Shi}.

In 1994 {\bf S-H Zhu} carefully went through Schoen-Yau's proof and showed
that in fact  $M^3$ with
$Ricci \ge 0$ that has $Ricci >0$ at one point
is diffeomorphic to $\RR^3$.  So the open case
in three dimensions was reduced to studying $M^3$ with a global vector field,
$V$, such that $Ricci(V,V) \ge 0$.  \cite {Zhu2}.  

Earlier {\bf S-H Zhu} had proven that if $M^3$ has only $Ricci \ge 0$ and 
\be
\lim_{r\to\infty}Vol(B_p(r))/r^3>0
\ee
then $M^3$ {\em is contractible} \cite{Zhu1}.  He first shows $\pi_1(M^3)=0$
because if $M^3=S^2 \times \RR$ then its volume growth is linear.
He then proves $\pi_1(M^3)$ is torsion free, which is true in general,
but contradicts Anderson and Li's assertion that it is finite unless 
$\pi_1(M^3)=0$ as well.  Zhu's Theorem is not true in dimensions four 
as demonstrated by Menguy and Otsu's examples \cite{Mng1} \cite{Ot}.

Without any volume assumptions,
we know from Shen-Sormani that $H_2(M^3,Z)$ is classified
and that $H_1(M,Z)$ is torsion free.  In fact $\pi_1(M)$
is completely understood as long as it is finitely generated.
From Wilking's reduction of the Milnor Conjecture, we then
need only  understand the topology of a three manifold, $N^3$, 
$\pi_1(M^3)$ infinitely generated and abelian.

Topologists {\bf Evans} and {\bf Moser} proved that if the fundamental group
of a three manifold
 is not finitely generated, then it is a subgroup of the
additive group of rational numbers \cite{EvMo}.  So we need only
worry about fundamental groups like the dyadic rationals:
\be \label{dyadiceqn}
\{ \, \frac{ p}{2^j} \, : \, p \in \ZZ, \,  j \in \{0,1,2,3,...\}\,\}.
\ee
It is important to note that there is a well-known topological 
construction of an open three manifold, $N^3$, whose 
fundamental group is the rationals.  It is based on the following example
studied thoroughly by {\bf Steenrod} in \cite{Str}
who credits {\bf Vietoris} with the idea.  A similar construction
by {\bf Whitehead} was used to construct a contractible three manifold
which isn't diffeomorphic to $\RR^3$ \cite{Wth}.

{\bf The Dyadic Solenoid Complement}
{\em
is an open topological manifold, $N^3$, such that
$
\pi_1(N^3)=\{k/2^j : \, k\in \ZZ, \, j\in \NN \}.
$
 Often
this example is described as $S^3$ with the dyadic solenoid
knot removed where the knot is a Cantor set bundle over an
$S^1$, but the following construction using embedded tori
gives a more geometric illustration.}

\begin{figure}[htbp]
\begin{center}
\includegraphics[height=2in ]{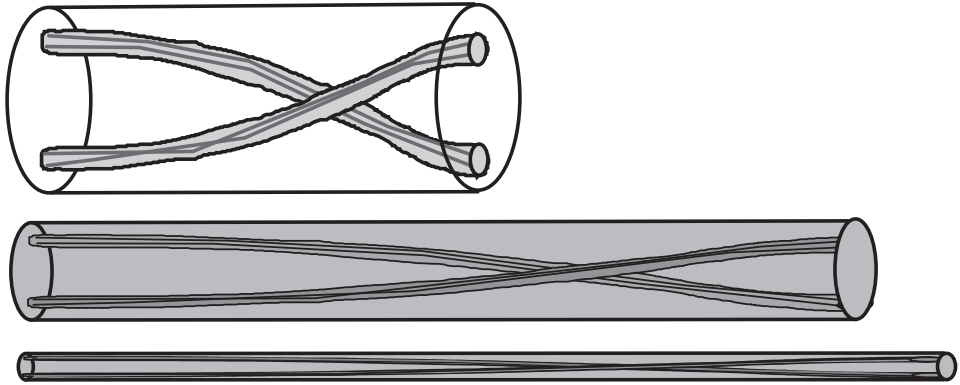}
\end{center}
\caption{Constructing the Dyadic Solenoid Complement} \label{3manifold}
\end{figure}

 {\bf Construction:}
We begin with a solid torus, $N_0=S^1 \times D^2$, with a curve 
$C_0:S^1 \to S^1 \times D^2$ with winding number 2.  Removing
a small neighborhood around the image of $C_0$, we get the building
block of our space
\be
N_i = (S^1\times D^2) \setminus T_{\epsilon} (C_0(S^1)).
\ee 
This block has two boundary components, $\partial N_i^-$ and
\be
\partial N_i^+= \partial T_\epsilon(C_0(S^1))
\ee
both of which are homeomorphic to $T^2$.    We build $N^3$ by
gluing all the $N_i$ together so that $\partial N_0$ is
glued to $\partial N_1^-$, and each $\partial N_i^+$ is glued
to $\partial N_{i+1}^-$.  This is depicted somewhat
inside-out in Figure~\ref{3manifold},
where we see each block $N_i$ is stretched out, wound around twice
and embedded into the previous block $N_{i-1}$.  

Note that the closure $K_j = Cl(\bigcup_{i=1}^j N_i)$ is compact and
that $K_j$ exhaust $N^3$.  Any curve $C:S^1 \to N^3$ is sitting
in some $K_j$.  Wrapping twice around $C$, one gets a curve homotopic
to a simple curve in $N_{j+1}\subset K_{j+1}$.  
In fact the halfway generators will have the
form: 
\be
g_1^2=g_0, \,\, g_2^2=g_1, \,\, g_3^2=g_2, \,\, g_4^2=g_3,...
\ee
if we choose a point $p \in N_0$ and make each subsequent $N_i$
significantly larger than the previous one.
The fundamental group of the dyadic solenoid is thus the dyadic rationals
with addition as the operation.

By taking the winding numbers of successive tori to run through every natural
number instead of just the number $2$, one obtains a space $N^3$ such
that $\pi_1(N^3)$ is the rationals.

{\bf Open Problem on the Dyadic Solenoid Complement:}
{\em By Schoen-Yau, we know $N^3$ does not admit a metric with positive
Ricci curvature but it is an open question as to whether it
admits a metric with $Ricci \ge 0$.}  Note that one must check it
satisfies all the topological restrictions on three manifolds with
$Ricci \ge 0$.

\sect{Open Problems:} \label{open}

We have attempted to state all known topological results and examples
concerning complete noncompact $M^n$ with $Ricci \ge 0$ 
or $Ricci >0$
that have either no additional conditions or only restrictions on:

* dimension

* diameter growth

* volume growth 

\noindent
There are a number of beautiful theorems with
additional conditions like bounds on conjugacy radius,
bounds on sectional curvature, quadratically
decaying positive lower bounds on Ricci curvature, and properness
to name a few.  However, without additional conditions, 
everything else is open.

For those wishing to investigate restrictions
to the topology of {\bf open four manifolds} with $Ricci \ge 0$
work by Noronha \cite{Nor} and Sha-Yang \cite{ShaYng3}
on compact four manifolds with $Ricci\ge 0$ might be helpful.  
Keep in mind the many examples given in this paper
as well as the fact that $H_3(M^4,Z)=0$ and $H_2(M^4,Z)$ is torsion
free except for $M^4$ with split double covers.

For those wishing to investigate restrictions of
the topology of manifolds with additional volume constraints,
one approach would be to assume $M^n$
has at most {\bf quadratic volume growth} and try to 
use Cheng-Yau's result that such an $M^n$ is parabolic 
to restrict the topology \cite{ChgYau}.

For those wishing to construct examples with interesting topology
remember
{\bf the key qualitative properties of $M^n$ are:}

* have a cover with an abelian fundamental group
\cite{Wlk}

* satisfy the loops to infinity property \cite{Sor4}

* have volume growth which is not maximal \cite{Li}
nor minimal \cite{Sor3}

* have large linear diameter growth \cite{Sor3}

\noindent
It is possible that the dyadic solenoid complement, $N^3$
satisfies all these conditions.

{\bf The Warped Dyadic Solenoid}, $N^3 \times \RR^k$,
is a more likely Milnor counter example than $N^3$ itself,
as we have seen how one can construct metrics
with $Ricci>0$ on vector bundles.  
Proving $N^3 \times \RR^k$ does not admit
a metric with $Ricci \ge 0$ would be of some interest
and might be done using harmonic map techniques developed by 
Schoen-Yau in \cite{SchYau2}.  Constructing a metric with 
$Ricci \ge 0$ on $N^3 \times \RR^k$
would be an astoundingly important result! It would lead to
many  new questions such as:

* What is the smallest diameter growth of a Milnor counter example?

* What is the largest volume growth of a Milnor counter example?

* What is the smallest volume growth of a Milnor counter example?

* What is the smallest dimension of a Milnor counter example?

\noindent
Even without the existence of a Milnor counter cxample, one may
well ask these questions of $M^n$ with a prescribed
almost nilpotent fundamental group.  

We hope that this survey article will provide new intuition
and insight allowing the readers to find open problems of
their own.  Until the topology of $M^n$ is completely understood, 
there is much work to be done.

\sect{Acknowledgments:}

We thank Guofang Wei and John Lott for clarifying some of the literature,
William Minicozzi for suggesting the quadratic volume growth problem,
Shing-Tung Yau for first mentioning to the second author that
$\pi_1(M^n)$ might be the rationals and Mikhail Khovanov for describing
the construction of the dyadic solenoid complement using embedded tori.
We would also  like to thank John Hempel, 
Bruce Kleiner, Rob Schneiderman, John Smillie and Rick Schoen for discussions
and emails related to the section on three manifolds.  The second author would
like to thank Courant Institute for its hospitality.

\end{document}